\RequirePackage{fix-cm}
\documentclass[smallextended]{svjour3}       % onecolumn (second format)
\smartqed  % flush right qed marks, e.g. at end of proof
\usepackage{graphicx}
\usepackage{epstopdf}
\usepackage{amssymb}
\usepackage{stfloats}
\usepackage{amsmath}
\usepackage{url}
\usepackage[colorlinks=true,citecolor=blue,linkcolor=blue]{hyperref}
\usepackage{multirow}
\usepackage{booktabs}
\begin{document}

\title{On a collocation method for the time-fractional convection-diffusion equation with variable coefficients}

\titlerunning{A collocation method for fractional convection-diffusion equation}        % if too long for running head

\author{X. G. Zhu \and Y. F. Nie}

%\authorrunning{Numer Algor} % if too long for running head

\institute{  X. G. Zhu \and Y. F. Nie\at
               Department of Applied Mathematics, Northwestern  Polytechnical University, Xi'an, Shanxi 710129, P. R. China \\   %%School of Natural and Applied Sciences
              \email{yfnie@nwpu.edu.cn}       %  \\
              \and
              X. G. Zhu \at
              \email{zhuxg590@yeah.net}           %  \\
}

\date{Received: date / Accepted: date}

\maketitle

\begin{abstract}
In this work, a new collocation approach using a combination of a wavelet operational matrix method and the exponential spline interpolation
is proposed to solve the time-fractional convection-diffusion equation with variable coefficients.
The operational matrix of fractional order integration is first derived based on sine-cosine wavelet functions,
which helps to convert the underlying equation into a linear algebraic system. Then, an exponential B-spline method is introduced in spatial direction.
On selecting a set of proper collocation points, the method in presence is evaluated on several test problems and the numerical results finally %at last
illustrate its validity and applicability. %robustness.
\end{abstract}

\section{Introduction}\label{intro}
As a generalization of classic models, fractional differential equations are encountered in a broad range of % a series of
individual disciplines recently, covering astrophysics, biology, electrochemistry, continuum mechanics, viscoelastic flows
%biology, hydrology,acoustics and psychology, continuum mechanics, finance, admit, fluid mechanics, reflect  standard ones   porous media
%fills the gap differential, global correlation, modeled in mathematical physics, shortcomings  chemistry
and so forth \cite{ref03,ref02,ref01}. These equations admit the non-local distributed effects and historical dependence, % of the system function,
which well remedy the serious drawback that the integer models can not apply to some of the natural phenomena like anomalous diffusion.
Due to the lack of analytic techniques and the closed form of true solutions,
numerical algorithms are placed  heavy reliance %is sent to great expectations. instrument great
and have gradually emerged as an essential tool utilized to study them \cite{Re05,ref08,ref04,ref07,ref06,ref05}.

In this regard, our aim is to derive an efficient method for solving the following convection-diffusion equation
\begin{align}\label{eq01}
    \frac{\partial^\alpha y(x,t)}{\partial t^\alpha}+a(x)\frac{\partial y(x,t)}{\partial x}
     +b(x)\frac{\partial^2 y(x,t)}{\partial x^2}=f(x,t),  \quad 0\leq x\leq \ell, \ 0<t\leq 1,
\end{align}
subjected to the initial and boundary conditions
\begin{align}
  & y(x,0)=\varphi(x),\quad 0\leq x\leq\ell, \label{eq02}\\
  & y(0,t)=g_1(t), \quad y(\ell,t)=g_2(t), \quad 0<t\leq 1, \label{eq03}
\end{align}
where $0<\alpha<1$, $a(x)$, $b(x)$ do not equal to zero at the same time, and $\varphi(x)$, $g_1(t)$, $g_2(t)$, $f(x,t)$ are properly %appropriately
prescribed. In Eq. (\ref{eq01}), the time-fractional derivative is defined in Caputo sense as %with enough smoothness.
\begin{align*}
 \frac{\partial^\alpha y(x,t)}{\partial t^\alpha}=\frac{1}{\Gamma(1-\alpha)}
  \int^t_0\frac{\partial y(x,\xi)}{\partial \xi}\frac{d\xi}{(t-\xi)^\alpha}, %J^{1-\alpha}\bigg(\frac{\partial y(x,t)}{\partial t}\bigg)=
\end{align*}
or alternatively by $J^{(1-\alpha)}\frac{\partial y(x,t)}{\partial t}$, where
\begin{align}\label{eq04}
  J^{(\mu)}y(x,t)=\frac{1}{\Gamma(\mu)}\int^t_0\frac{y(x,\xi)d\xi}{(t-\xi)^{1-\mu}},
\end{align}
which defines the $\mu$-th Riemann-Liouville integration and coincides with the $n$-th classic integration of the form %enjoys
\begin{align*}
  J^{(n)}y(x,t)&=\int^t_0d\xi_1\int^{\xi_1}_0d\xi_2 \cdots \int^{\xi_{n-1}}_0y(x,\xi_n)d\xi_n,
 %   &=\frac{1}{(n-1)!}\int^t_0 (t-\xi)^{n-1}y(x,\xi)d\xi,
\end{align*}
when $\mu=n\in \mathbb{N}$. In particular, we note that Eq. (\ref{eq04}) recovers to the common used definite
integration with the integral range $[0,t]$ if $\mu=1$ is taken.  %recovers degenerates

Apart from a few analytic approximations, many robust numerical algorithms have been designed
and applied to solve Eqs. (\ref{eq01})-(\ref{eq03}). In \cite{ref14,ref13,ref10,ref12,ref11,ref09},
the authors considered this problem of diffusion type with constant coefficients by using the implicit difference method,
Legendre spectral method, high-order compact difference method, finite element method,
direct discontinuous Galerkin method, and the moving least squares implicit meshless method, respectively.  % radial basis
Chen et al.\ proposed a Haar wavelet operational matrix method for Eqs. (\ref{eq01})-(\ref{eq03}) \cite{ref15}.
Uddin and Haq used a radial basis interpolation method to deal with this equation with constant coefficients \cite{ref16}.
In \cite{ref17}, a collocation technique using Sinc functions and shifted Legendre polynomials was established
and a spectral method for Gegenbauer polynomials can be found in \cite{ref18}.
In \cite{ref21,ref20}, the operational matrix methods based on 2D-Block Pulse and shifted Jacobi polynomial functions
were proposed for the fractional partial differential equations. % with constant or variable coefficients.
Luo et al.\  developed a quadratic spline collocation method for the constant coefficient case without convection \cite{ref22}.
Sayevand et al.\ considered the same type equation by a cubic B-spline collocation method \cite{ref23}.
Pirkhedri and Javadi gave a collocation approach for the variable coefficient case
via expanding its solution as the elements of Haar and Sinc functions \cite{ref24}.
In this article, inspired by the works above, we tackle Eqs. (\ref{eq01})-(\ref{eq03}) on the basis of
sine-cosine wavelets and exponential spline trial functions via a collocation strategy.
The wavelet operational matrix is derived and with it the weakly singular fractional integration is
eliminated by reducing the original problem to those of solving a set of algebraic equations.
The stated method calls for a low time cost and is apt to be realized by computer.

The layout is organized as follows. In Section \ref{Se1}, we describe the sine-cosine wavelet functions, %orthogonal
exponential B-spline basis, and some of their basic properties. The sine-cosine wavelets operational matrix of fractional integration
and the proposed collocation based method are investigated in Section \ref{Se2} and Section \ref{Se3}. Illustrative examples
are included in Section \ref{Se4}. The last section outlines conclusions.

\section{Preliminaries} \label{Se1}
The following basic definitions will be very useful hereinafter.
\subsection{Wavelets and sine-cosine wavelets}
Wavelets are a sequence of rescaled functions constructed from dilation and translation of a scaling function called
the mother wavelet. If the dilation and translation factors $a$, $b$ are continuous, we have
the continuous wavelets \cite{ref25}:
\begin{align*}
\psi_{a,b}(t)=|a|^{-\frac{1}{2}}\psi\bigg(\frac{t-b}{a}\bigg), \quad a,b \in \mathbb{R}, \quad a\neq0,
\end{align*}
whereas we have the discrete wavelets:
\begin{align*}
\psi_{k,n}(t)=|a_0|^{\frac{k}{2}}\psi(a_0^{k}t-nb_0),
\end{align*}
if the factors $a$, $b$ are limited to discrete values as $a=a_0^{-k}$, $b=nb_0a_0^{-k}$, %basis restricted
where $a_0>1$, $b_0>0$, and $n$, $k$ are positive integers. The discrete wavelets form a compact wavelet basis on $L^2(\mathbb{R})$ and
contain an orthonormal basis as a special case when $a_0=2$ and $b_0=1$.

The orthonormal sine-cosine wavelets $\psi_{n,m}(t)$ are defined on $[0,1)$ as follow
%Sine-cosine wavelets $\psi_{n,m}(t)$ are orthonormal on $[0,1)$ and defined by
\begin{align*}
\psi_{n,m}(t)=\left\{
\begin{array}{lc}
2^{\frac{k+1}{2}}\textrm{scw}_m(2^kt-n),&\quad\frac{n}{2^k}\leq t<\frac{n+1}{2^k},\\
0,& \textrm{otherwise},
\end{array} \right.
\end{align*}
with
\begin{align*}
\qquad\quad\textrm{scw}_m(t)=\left\{
\begin{array}{ll}
\frac{1}{\sqrt{2}},&\quad m=0,\\
\textrm{cos}(2m\pi t),&\quad m=1,2,\ldots,M,\\   %1\leq m\leq L,\\
\textrm{sin}(2(m-M)\pi t),&\quad m=M+1,M+2,\ldots,2M,   %L+1\leq m\leq 2L,\\
\end{array} \right.
\end{align*}
where $n=0,1,\ldots,2^k-1$, $k\in \mathbb{N}\cup \{0\}$, $M\in \mathbb{Z}$ and $t$ is the normalized time.

\subsection{Exponential spline functions}
Given a set of equidistant knots $\{x_j\}_{j=0}^{N_h}$ on $[0,\ell]$ together with another six ghost
knots $x_j$, $j=-3,-2,-1,N_h+1,N_h+2,N_h+3$ beyond $[0,\ell]$, we denote
\begin{align*}
 h=\ell/N_h, \quad  p=\max\limits_{0\leq j\leq N_h}p_j,
 \quad s=\sinh(ph), \quad c=\cosh(ph),
\end{align*}
where $N_h\in\mathbb{Z}^+$ and $p_j$ is the value of $p(x)$ at $x_j$.
The exponential splines are recognized as an extension of the semi-classic cubic splines, which
also admit a basis of exponential B-splines $\{B_j(x)\}_{j=-1}^{N_h+1}$ as follows \cite{ref26}

\begin{align*}
B_j(x)=\left\{
\begin{aligned}
&e(x_{j-2}-x)-\frac{e}{p}\sinh(p(x_{j-2}-x)), \, \qquad\quad\qquad\qquad\qquad x\in [x_{j-2},x_{j-1}], \\
&a+b(x_j-x)+c\exp(p(x_j-x))+d\exp(-p(x_j-x)), \quad  x\in [x_{j-1},x_{j}], \\
&a+b(x-x_j)+c\exp(p(x-x_j))+d\exp(-p(x-x_j)), \quad x\in [x_{j},x_{j+1}], \\
&e(x-x_{j+2})-\frac{e}{p}\sinh(p(x-x_{j+2})), \, \qquad\quad\qquad\qquad\qquad x\in [x_{j+1},x_{j+2}], \\
&0, \quad \textrm{otherwise},
\end{aligned}
\right.
\end{align*}
where
\begin{gather*}
e=\frac{p}{2(phc-s)},\quad a=\frac{phc}{phc-s},\quad b=\frac{p}{2}\bigg[\frac{c(c-1)+s^2}{(phc-s)(1-c)}\bigg],\\
c=\frac{1}{4}\bigg[\frac{\exp(-ph)(1-c)+s(\exp(-ph)-1)}{(phc-s)(1-c)}\bigg], \quad
d=\frac{1}{4}\bigg[\frac{\exp(ph)(c-1)+s(\exp(ph)-1)}{(phc-s)(1-c)}\bigg].
\end{gather*}
The values of $B_j(x)$ at each knot are given as
\begin{align}\label{eq23}
B_j(x_k)=\left\{
\begin{aligned}
&1,  &\textrm{if}\ \ k=j, \\
&\frac{s-ph}{2(phc-s)},  &\textrm{if} \ \ k=j\pm1,\\
&0,  &\textrm{if} \ \ k=j\pm2.
\end{aligned}
\right.
\end{align}
The values of $B'_j(x)$ and $B''_j(x)$ at each knot are given as
\begin{align}\label{eq24}
B'_j(x_k)=\left\{
\begin{aligned}
&0  &\textrm{if}\ \ k=j, \\
&\frac{\mp p(1-c)}{2(phc-s)},  &\textrm{if} \ \ k=j\pm1,\\
&0,  &\textrm{if} \ \ k=j\pm2,
\end{aligned}
\right.
\end{align}
and
\begin{align}\label{eq25}
B''_j(x_k)=\left\{
\begin{aligned}
&\frac{-p^2s}{phc-s},  &\textrm{if}\ \ k=j, \\
&\frac{p^2s}{2(phc-s)},  &\textrm{if} \ \ k=j\pm1,\\
&0,  &\textrm{if} \ \ k=j\pm2.
\end{aligned}
\right.
\end{align}
It is noteworthy that $\{B_j(x)\}_{j=-1}^{N_h+1}$ form an exponential spline space on $[0,\ell]$,
each of which is twice continuously differentiable in respect of $x$.

\section{Operational matrix of fractional integration}\label{Se2}
In this part, let $N_t=2^k(2M+1)$ and
\begin{align*}
\mathbf{\Psi}(t)=[\psi_{0,0},\psi_{0,1},\ldots,\psi_{0,2M},\psi_{1,0},\ldots,\psi_{1,2M},\ldots,\psi_{2^k-1,2M}],
\end{align*}
we shall derive a matrix $\textbf{J}^\mu_{N_t\times N_t}$ such that %shall
%\begin{align}\label{eq05}
%J^{(\mu)}\mathbf{\Psi}^T(t)=\frac{1}{\Gamma(\mu)}\int^t_0\frac{\mathbf{\Psi}^T(\xi)d\xi}{(t-\xi)^{1-\mu}}
%\cong\textbf{J}^\mu_{N_t\times N_t}\mathbf{\Psi}^T(t),
%\end{align}
\begin{align}\label{eq05}
J^{(\mu)}\mathbf{\Psi}(t)=\frac{1}{\Gamma(\mu)}\int^t_0\frac{\mathbf{\Psi}(\xi)d\xi}{(t-\xi)^{1-\mu}}
\cong\mathbf{\Psi}(t)\textbf{J}^{\mu,T}_{N_t\times N_t},
\end{align}
%which is termed the operational matrix of fractional integration. To begin with,
%we recall the $N_t$-set of Block-Pulse functions:
%\begin{align*}
%b_i(t)=\left\{
%\begin{array}{ll}
%1,&\quad \frac{i}{N_t}\leq t<\frac{i+1}{N_t},\\
%0,&\quad \textrm{otherwise},
%\end{array} \right.
%\end{align*}
%with $i=0,1,\ldots,N_t-1$.
%They are disjoint and orthogonal over $[0,1)$.
termed the operational matrix of integration $J^{(\mu)}$. To begin with,
we recall the $N_t$-set of Block-Pulse functions $b_i(t)$, $i=0,1,\ldots,N_t-1$, which equal to $1$ if $\frac{i}{N_t}\leq t<\frac{i+1}{N_t}$
and vanish otherwise. They are disjoint and orthogonal over $[0,1)$.
For each $\psi_{n,m}(t)$, we now approximate it by the truncated series
\begin{align*}
  \psi_{n,m}(t)\cong\sum_{i=0}^{N_t-1}Q_ib_i(t), \quad Q_i=N_t\int_0^1b_i(\xi)\psi_{n,m}(\xi)d\xi.
\end{align*}
Using the compactness of $b_i(t)$, we have
\begin{align*}
 Q_i=N_t\int_{i/N_t}^{(i+1)/N_t}b_i(\xi)\psi_{n,m}(\xi)d\xi\cong
  \psi_{n,m}\bigg(\frac{2i+1}{2N_t}\bigg),
\end{align*}
Then, a relationship between the sine-cosine wavelets and Block-Pulse functions can be built, i.e.,
\begin{align}\label{eq06}
 \mathbf{\Psi}^T(t)\cong\textbf{Q}\textbf{B}(t),
\end{align}
where $\textbf{Q}=\textrm{diag}(\textbf{Q}_0,\textbf{Q}_1,\ldots,\textbf{Q}_{2^k-1})$, $\textbf{B}(t)=[b_0(t),b_1(t),\ldots,b_{N_t-1}(t)]^T$, and
\begin{align*}
     \textbf{Q}_n=\left(
    \begin{array}{ccccc}
      \psi_{n,0}\big(\frac{\tilde{n}+1}{2N_t}\big) & \psi_{n,0}\big(\frac{\tilde{n}+3}{2N_t}\big)& \cdots & \psi_{n,0}\big(\frac{\tilde{n}+4M+1}{2N_t}\big)\\
      \psi_{n,1}\big(\frac{\tilde{n}+1}{2N_t}\big) & \psi_{n,1}\big(\frac{\tilde{n}+3}{2N_t}\big)& \cdots & \psi_{n,1}\big(\frac{\tilde{n}+4M+1}{2N_t}\big)\\
      \psi_{n,2}\big(\frac{\tilde{n}+1}{2N_t}\big) & \psi_{n,2}\big(\frac{\tilde{n}+3}{2N_t}\big)& \cdots & \psi_{n,2}\big(\frac{\tilde{n}+4M+1}{2N_t}\big)\\
      \vdots& \vdots & & \vdots&  \\
      \psi_{n,2M}\big(\frac{\tilde{n}+1}{2N_t}\big) & \psi_{n,2M}\big(\frac{\tilde{n}+3}{2N_t}\big)&\cdots & \psi_{n,2M}\big(\frac{\tilde{n}+4M+1}{2N_t}\big)\\
    \end{array}
  \right),
\end{align*}
with $n=0,1,\ldots,2^k-1$ and $\tilde{n}=2n(2M+1)$. $\textbf{Q}_n$ are a family of $(2M+1)\times(2M+1)$-dimensional matrices.

The Block-Pulse fractional operational matrix is as follow
\begin{equation}\label{eq22}
J^{(\mu)}\mathbf{B}(t)\cong\textbf{F}^\mu_{N_t\times N_t}\mathbf{B}(t),
\end{equation}
provided by Kilicman and Al Zhour  \cite{ref19}, where
\begin{align*}
\textbf{F}^\mu_{N_t\times N_t}=\frac{1}{N_t^\mu}\frac{1}{\Gamma(\mu+2)}\left(
    \begin{array}{cccccc}
      1 &\xi_1 &\xi_2 &\xi_3 &\ldots &\xi_{N_t-1} \\
      0 &1 &\xi_1 &\xi_2 &\ldots &\xi_{N_t-2} \\
      0 &0 &1 &\xi_1 &\ldots &\xi_{N_t-3} \\
      \vdots & &\ddots & &\ddots &\vdots\\
      0 &0 &\ldots &0 &1 &\xi_1\\
      0 &0 &0 &\ldots &0 &1
    \end{array}
  \right),
\end{align*}
and $\xi_k=(k+1)^{\mu+1}-2k^{\mu+1}+(k-1)^{\mu+1}$. Then, we obtain
\begin{equation*}
    J^{(\mu)}\mathbf{\Psi}^T(t)\cong\textbf{Q}J^{(\mu)}\textbf{B}(t)
    \cong\textbf{Q}\textbf{F}^\mu_{N_t\times N_t}\textbf{B}(t)\cong\textbf{J}^\mu_{N_t\times N_t}\mathbf{\Psi}^T(t),
\end{equation*}
from Eqs. (\ref{eq05})-(\ref{eq22}), which implies
$\textbf{J}^\mu_{N_t\times N_t}=\textbf{Q}\textbf{F}^\mu_{N_t\times N_t}\textbf{Q}^{-1}$ finally. In particular,
we have for $\mu=0.5$, $k=1$, and $M=1$, that
\begin{align*}
\textbf{J}^{0.5}_{6\times6}=\left(
    \begin{array}{cccccc}
    0.5319   &-0.0209   &-0.1715  &  0.4407   & 0.0180  &  0.0821\\
   -0.0209   & 0.1651   & 0.0991  &  0.0180   & 0.0061  &  0.0148\\
    0.1715   &-0.0991   & 0.2243  & -0.0821   &-0.0148  & -0.0449\\
         0   &      0   &      0  &  0.5319   &-0.0209  & -0.1715\\
         0   &      0   &      0  & -0.0209   & 0.1651  &  0.0991\\
         0   &      0   &      0  &  0.1715   &-0.0991  &  0.2243
    \end{array}
  \right).
\end{align*}

\section{Description of the proposed method}\label{Se3}
In the sequel, we propose a matrix method for Eqs. (\ref{eq01})-(\ref{eq03}) by utilizing $\textbf{J}^\mu_{N_t\times N_t}$
and the exponential B-spline basis via a collocation technique.

\subsection{Function approximation}
Define $V_{N_h+3}=\textrm{span}\{B_{-1}(x),B_0(x),\ldots,B_{N_h+1}(x)\}$
as a $(N_h+3)$-dimensional exponential spline space on a uniform spatial partition with %$h=\ell/N_h$ and
the knots $\{x_j\}_{j=-3}^{N_h+3}$ inside or outside $[0,\ell]$ as before. If $y(x,t)$ is the exact solution of
Eqs. (\ref{eq01})-(\ref{eq03}), then we can expand $\frac{\partial y(x,t)}{\partial t}$ into a finite series as %its derivative
\begin{equation}\label{eq07}
\frac{\partial y_{N}(x,t)}{\partial t}=\sum_{n=0}^{2^k-1}\sum_{m=0}^{2M}\sum_{l=-1}^{N_h+1}c_{n,m,l}\psi_{n,m}(t)B_l(x).
\end{equation}
Integrating Eq. (\ref{eq07}) from $0$ to $t$, an approximate solution in the form
\begin{equation}\label{eq08}
y_N(x,t)=\sum_{n=0}^{2^k-1}\sum_{m=0}^{2M}\sum_{l=-1}^{N_h+1}c_{n,m,l}J^{(1)}\psi_{n,m}(t)B_l(x)+\varphi(x),
\end{equation}
to the model problem is obtained and is sought on $L^2[0,1]\cap V_{N_h+3}$ with the unknown weights $c_{n,m,l}$ yet to be determined,
where $\varphi(x)$ is the initial function. Along the same line, we can get %fashion
\begin{align}
\frac{\partial^\alpha y_N(x,t)}{\partial t^\alpha}&=\sum_{n=0}^{2^k-1}\sum_{m=0}^{2M}
    \sum_{l=-1}^{N_h+1}c_{n,m,l}J^{(1-\alpha)}\psi_{n,m}(t)B_l(x),\label{eq09}\\
\frac{\partial y_N(x,t)}{\partial x}&=\sum_{n=0}^{2^k-1}\sum_{m=0}^{2M}\sum_{l=-1}^{N_h+1}c_{n,m,l}J^{(1)}\psi_{n,m}(t)
    \frac{\partial B_l(x)}{\partial x}+\frac{\partial\varphi(x)}{\partial x},\label{eq10}\\
\frac{\partial^2 y_N(x,t)}{\partial x^2}&=\sum_{n=0}^{2^k-1}\sum_{m=0}^{2M}\sum_{l=-1}^{N_h+1}c_{n,m,l}J^{(1)}\psi_{n,m}(t)
    \frac{\partial^2 B_l(x)}{\partial x^2}+\frac{\partial^2 \varphi(x)}{\partial x^2}, \label{eq11}
\end{align}
on acting the operators $J^{(1-\alpha)}$, $\frac{\partial}{\partial x}$, and $\frac{\partial^2}{\partial x^2}$
on both sides of Eq. (\ref{eq07}) or Eq. (\ref{eq08}), respectively. Let
\begin{align*}
&\textbf{C}=[\textbf{C}_{-1},\textbf{C}_0,\ldots,\textbf{C}_{N_h+1}]^T,\\
&\textbf{C}_l=[c_{0,0,l},c_{0,1,l},\ldots,c_{0,2M,l},c_{1,0,l},\ldots,c_{1,2M,l},\ldots,c_{2^k-1,2M,l}],\\
&\textbf{H}(x)=[B_{-1}(x),B_{0}(x),\ldots,B_{N_h+1}(x)],\\
\end{align*}
with $l=0,1,\ldots,N_h$. Then Eqs. (\ref{eq08})-(\ref{eq11}) are reduced to
\begin{gather}
y_N(x,t)=\textbf{H}(x)\otimes\mathbf{\Psi}(t)\textbf{J}^{1,T}_{N_t\times N_t}\cdot\textbf{C}+\varphi(x),\label{eq12}\\
\frac{\partial^\alpha y_N(x,t)}{\partial t^\alpha}=\textbf{H}(x)\otimes\mathbf{\Psi}(t)\textbf{J}^{1-\alpha,T}_{N_t\times N_t}
    \cdot\textbf{C},\label{eq13}\\
\frac{\partial y_N(x,t)}{\partial x}=\frac{\partial\textbf{H}(x)}{\partial x}\otimes\mathbf{\Psi}(t)\textbf{J}^{1,T}_{N_t\times N_t}
    \cdot\textbf{C} + \frac{\partial \varphi(x)}{\partial x},\label{eq14}\\
\frac{\partial^2 y_N(x,t)}{\partial x^2}=\frac{\partial^2\textbf{H}(x)}{\partial x^2}\otimes\mathbf{\Psi}(t)\textbf{J}^{1,T}_{N_t\times N_t}
    \cdot\textbf{C} + \frac{\partial^2 \varphi(x)}{\partial x^2},\label{eq15}
\end{gather}
%\begin{align}
%&y_N(x,t)=\textbf{H}(x)\otimes\mathbf{\Psi}(t)\textbf{J}^{1,T}_{N_t\times N_t}\cdot\textbf{C}+\varphi(x),\label{eq12}\\
%&\frac{\partial^\alpha y_N(x,t)}{\partial t^\alpha}=\textbf{H}(x)\otimes\mathbf{\Psi}(t)\textbf{J}^{1-\alpha,T}_{N_t\times N_t}
%    \cdot\textbf{C},\label{eq13}\\
%&\frac{\partial y_N(x,t)}{\partial x}=\frac{\partial\textbf{H}(x)}{\partial x}\otimes\mathbf{\Psi}(t)\textbf{J}^{1,T}_{N_t\times N_t}
%    \cdot\textbf{C} + \frac{\partial \varphi(x)}{\partial x},\label{eq14}\\
%&\frac{\partial^2 y_N(x,t)}{\partial x^2}=\frac{\partial^2\textbf{H}(x)}{\partial x^2}\otimes\mathbf{\Psi}(t)\textbf{J}^{1,T}_{N_t\times N_t}
%    \cdot\textbf{C} + \frac{\partial^2 \varphi(x)}{\partial x^2},\label{eq15}
%\end{align}
in the matrix-vector forms thanks to the operational matrices $\textbf{J}^{1-\alpha}_{N_t\times N_t}$ and $\textbf{J}^{1}_{N_t\times N_t}$,
where the symbol ``$\otimes$'' denotes the Kronecker product.

\subsection{Construction of the collocation method}
On inserting the collocation points
\begin{gather}
\begin{aligned}\label{eq21}
  &t_i=\frac{2i-1}{2^k(2M+1)},\quad i=1,2,\ldots,N_t,\\
  & x_j=jh,\quad j=0,1,\ldots,N_h,
\end{aligned}
\end{gather}
and Eqs. (\ref{eq13})-(\ref{eq15}) into Eq. (\ref{eq01}), we have
\begin{align*}
&\textbf{H}(x_j)\otimes\mathbf{\Psi}(t_i)\textbf{J}^{1-\alpha,T}_{N_t\times N_t}\cdot\textbf{C}
+a(x_j)\frac{\partial\textbf{H}(x_j)}{\partial x}\otimes\mathbf{\Psi}(t_i)\textbf{J}^{1,T}_{N_t\times N_t}\cdot\textbf{C} \\
&\ +b(x_j)\frac{\partial^2\textbf{H}(x_j)}{\partial x^2}\otimes\mathbf{\Psi}(t_i)\textbf{J}^{1,T}_{N_t\times N_t}
\cdot\textbf{C}=f(x_j,t_i)-a(x_j)\frac{\partial \varphi(x_j)}{\partial x}-b(x_j)\frac{\partial^2 \varphi(x_j)}{\partial x^2},
\end{align*}
which admits a set of algebraic equations
\begin{align}\label{eq16}
\textbf{A}\otimes \textbf{P}^\alpha\textbf{C} +\textrm{diag}(\textbf{a})\textbf{K}_1\otimes \textbf{P}\textbf{C}
+\textrm{diag}(\textbf{b})\textbf{K}_2\otimes \textbf{P}\textbf{C}=\textbf{F}
\end{align}
due to the properties (\ref{eq23})-(\ref{eq25}) of exponential B-spline basis, where
\begin{equation*}
\textbf{a}=[a(x_0),a(x_1),\ldots,a(x_{N_h})],\ \ \textbf{b}=[b(x_0),b(x_1),\ldots,b(x_{N_h})],
\end{equation*}
\begin{align*}
\textbf{P}^\alpha=\left(
\begin{array}{c}
  \mathbf{\Psi}(t_1)\textbf{J}^{1-\alpha,T}_{N_t\times N_t} \\
  \mathbf{\Psi}(t_2)\textbf{J}^{1-\alpha,T}_{N_t\times N_t} \\
  \mathbf{\Psi}(t_3)\textbf{J}^{1-\alpha,T}_{N_t\times N_t} \\
  \vdots\\
  \mathbf{\Psi}(t_{N_t})\textbf{J}^{1-\alpha,T}_{N_t\times N_t}
\end{array}
\right), \ \
\textbf{P}=\left(
\begin{array}{c}
  \mathbf{\Psi}(t_1)\textbf{J}^{1,T}_{N_t\times N_t} \\
  \mathbf{\Psi}(t_2)\textbf{J}^{1,T}_{N_t\times N_t} \\
  \mathbf{\Psi}(t_3)\textbf{J}^{1,T}_{N_t\times N_t} \\
  \vdots\\
  \mathbf{\Psi}(t_{N_t})\textbf{J}^{1,T}_{N_t\times N_t}
\end{array}
\right),
\end{align*}
\begin{align*}
&\textbf{A}=\textrm{tri}\bigg[\frac{s-ph}{2(phc-s)},1,\frac{s-ph}{2(phc-s)}\bigg],\\
&\textbf{K}_1=\textrm{tri}\bigg[\frac{p(1-c)}{2(phc-s)},0,-\frac{p(1-c)}{2(phc-s)}\bigg],\\
&\textbf{K}_2=\textrm{tri}\bigg[\frac{p^2s}{2(phc-s)},-\frac{p^2s}{phc-s},\frac{p^2s}{2(phc-s)}\bigg],
\end{align*}
and the right-hand vector
\begin{gather*}
\textbf{F}=\textbf{f}-\textrm{diag}(\textbf{a})\boldsymbol{\varphi}_x\otimes \textbf{I}_{N_t\times1}
  -\textrm{diag}(\textbf{b})\boldsymbol{\varphi}_{xx}\otimes \textbf{I}_{N_t\times1},\\
\textbf{I}_{N_t\times1}=[1,\ 1,\ \ldots,\ 1]^T_{N_t\times1},\\
\boldsymbol{\varphi}_x=\bigg[\frac{\partial\varphi(x_0)}{\partial x},\frac{\partial\varphi(x_1)}{\partial x},\ldots,
    \frac{\partial\varphi(x_{N_h})}{\partial x}\bigg]^T,\\
\boldsymbol{\varphi}_{xx}=\bigg[\frac{\partial^2\varphi(x_0)}{\partial x^2},\frac{\partial^2\varphi(x_1)}{\partial x^2},
    \ldots,\frac{\partial^2\varphi(x_{N_h})}{\partial x^2}\bigg]^T,\\
\quad\textbf{f}=[f(x_0,t_1),f(x_0,t_2),\ldots,f(x_0,t_{N_t}),f(x_1,t_1),\ldots,f(x_{N_h},t_{N_t})]^T.
\end{gather*}
The ``tri" means generating a tri-diagonal matrix of size $(N_h+1)\times(N_h+3)$.

However, the system (\ref{eq16}) comprises $N_t(N_h+1)$ equations in of $N_t(N_h+3)$ unknowns, which is not solvable.
In order to avoid this issue, the boundary conditions Eq. (\ref{eq03}) is further utilized so that another %utilized
$2N_t$ additional constraints can be derived. Taking $x=x_0$ and $x=x_{N_h}$ in Eq. (\ref{eq12}) leads to
\begin{align}
&\textbf{H}(x_0)\otimes\mathbf{\Psi}(t)\textbf{J}^{1,T}_{N_t\times N_t}\cdot\textbf{C}+\varphi(x_0)=g_1(t),\label{eq17}\\
&\textbf{H}(x_{N_h})\otimes\mathbf{\Psi}(t)\textbf{J}^{1,T}_{N_t\times N_t}\cdot\textbf{C}+\varphi(x_{N_h})=g_2(t).\label{eq18}
\end{align}
Inserting the collocation points into Eqs. (\ref{eq17})-(\ref{eq18}), it turns out that
\begin{align}
\textbf{Z}_1\otimes\textbf{P}\textbf{C}&=\textbf{g}_{1}-\varphi(x_0)\textbf{I}_{N_t\times1},\label{eq19}\\
\textbf{Z}_2\otimes\textbf{P}\textbf{C}&=\textbf{g}_{2}-\varphi(x_{N_h})\textbf{I}_{N_t\times1},\label{eq20}
\end{align}
by making use of the properties (\ref{eq23})-(\ref{eq25}), where
\begin{align*}
&\textbf{g}_\nu=[g_\nu(t_1),g_\nu(t_2),\cdots,g_\nu(t_{N_t})]^T,\quad \nu=1,2,\\
&\textbf{Z}_1=\bigg[\frac{s-ph}{2(phc-s)},1,\frac{s-ph}{2(phc-s)},0,\ldots,0\bigg]_{1\times(N_h+3)},\\
&\textbf{Z}_2=\bigg[0,\ldots,0,\frac{s-ph}{2(phc-s)},1,\frac{s-ph}{2(phc-s)}\bigg]_{1\times(N_h+3)}.
\end{align*}
Recombining Eqs. (\ref{eq16}), (\ref{eq19})-(\ref{eq20}), $\textbf{C}$ can thus be computed and being substituted into
Eq. (\ref{eq12}), an approximate solution to Eqs. (\ref{eq01})-(\ref{eq03}) is finally obtained.

\section{Illustrative examples}\label{Se4}
In this section, numerical examples are performed to illustrate the actual performance of our proposed method.
In the tests, the errors are measured by
%\begin{align*}
%&||\textrm{error}||_{L^2}=\sqrt{h\sum^{N_h-1}_{j=1}\Big|y(x_j,t)-y_N(x_j,t)\Big|^2},\\
%&||\textrm{error}||_{L^\infty}=\max\limits_{1\leq j\leq N_h-1}\Big|y(x_j,t)-y_N(x_j,t)\Big|,
%\end{align*}
\begin{align*}
&e_2(t,N_h)=\sqrt{\frac{\ell}{N_h}\sum^{N_h-1}_{j=1}\Big|y(x_j,t)-y_N(x_j,t)\Big|^2},\\
&e_\infty(t,N_h)=\max\limits_{1\leq j\leq N_h-1}\Big|y(x_j,t)-y_N(x_j,t)\Big|,
\end{align*}
with a fixed $t$ and the convergent rates are computed by
\begin{align*}
\textrm{Cov. rate}=\textrm{log}_2\bigg(\frac{e_\nu(t,N_h)}{e_\nu(t,2N_h)}\bigg), \quad \nu=1,\infty.
\end{align*}
The free parameter $p$ is chosen and its optimal value for a concrete problem is numerically studied.
We implement the method as the following algorithm:
\begin{center}
\begin{tabular}{ll}
  \toprule
  \multicolumn{2}{l}{\textbf{Algorithm}} \\
  \midrule
  1. & Input $\alpha$, $p$, $k$, $M$, and $N_h$\\
  2. & Initialize \textbf{A}, $\textbf{K}_1$, $\textbf{K}_2$, $\textbf{P}^\alpha$, $\textbf{P}$, and $\textbf{F}$\\
  3. & Allocate the collocation points as (\ref{eq21})\\
  4. & Construct the operational matrices $\textbf{J}^{1-\alpha}_{N_t\times N_t}$ and $\textbf{J}^{1}_{N_t\times N_t}$\\
  5. & Compute $\textbf{A}\otimes \textbf{P}^\alpha$, $\textbf{K}_1\otimes \textbf{P}$,
               $\textbf{K}_2\otimes \textbf{P}$, and $\textbf{F}$ to get Eq. (\ref{eq16})\\
  6. & Compute $\textbf{Z}\otimes \textbf{P}^\alpha$, $\textbf{Z}_1\otimes \textbf{P}$, $\textbf{g}_1$,
               and $\textbf{g}_2$ to get Eqs. (\ref{eq19})-(\ref{eq20}) \\
  7. & Reform  Eqs. (\ref{eq16}) and (\ref{eq19})-(\ref{eq20}) and solve the system for $\textbf{C}$ \\
  8. & Output $y_N(x,t)=\textbf{H}(x)\otimes\mathbf{\Psi}(t)\textbf{J}^{1,T}_{N_t\times N_t}\cdot\textbf{C}+\varphi(x)$\\
  \bottomrule
\end{tabular}
\end{center}

\begin{figure}
\begin{minipage}[t]{0.32\linewidth}
\includegraphics[width=1.65in]{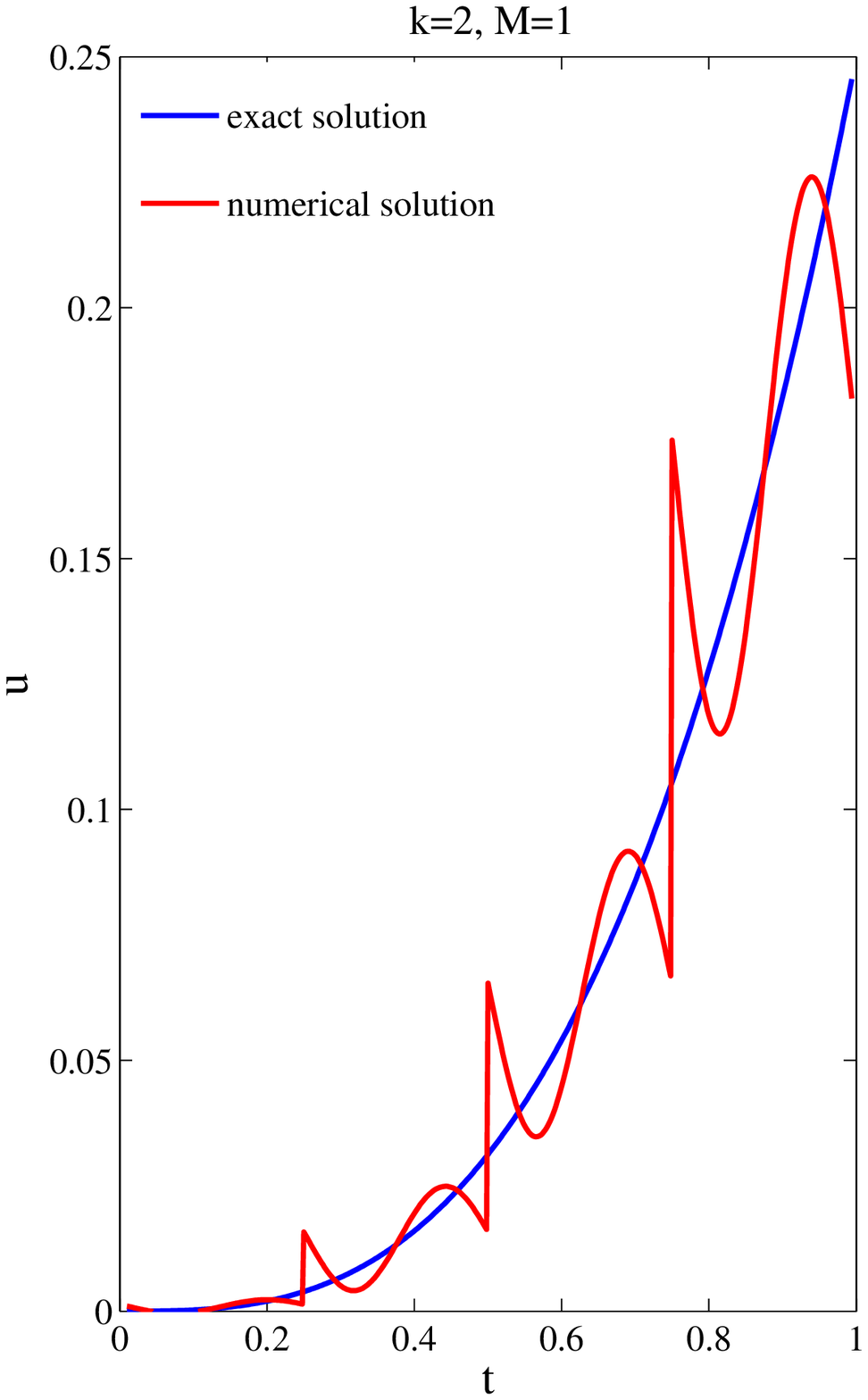}
\end{minipage}
\begin{minipage}[t]{0.32\linewidth}
\includegraphics[width=1.65in]{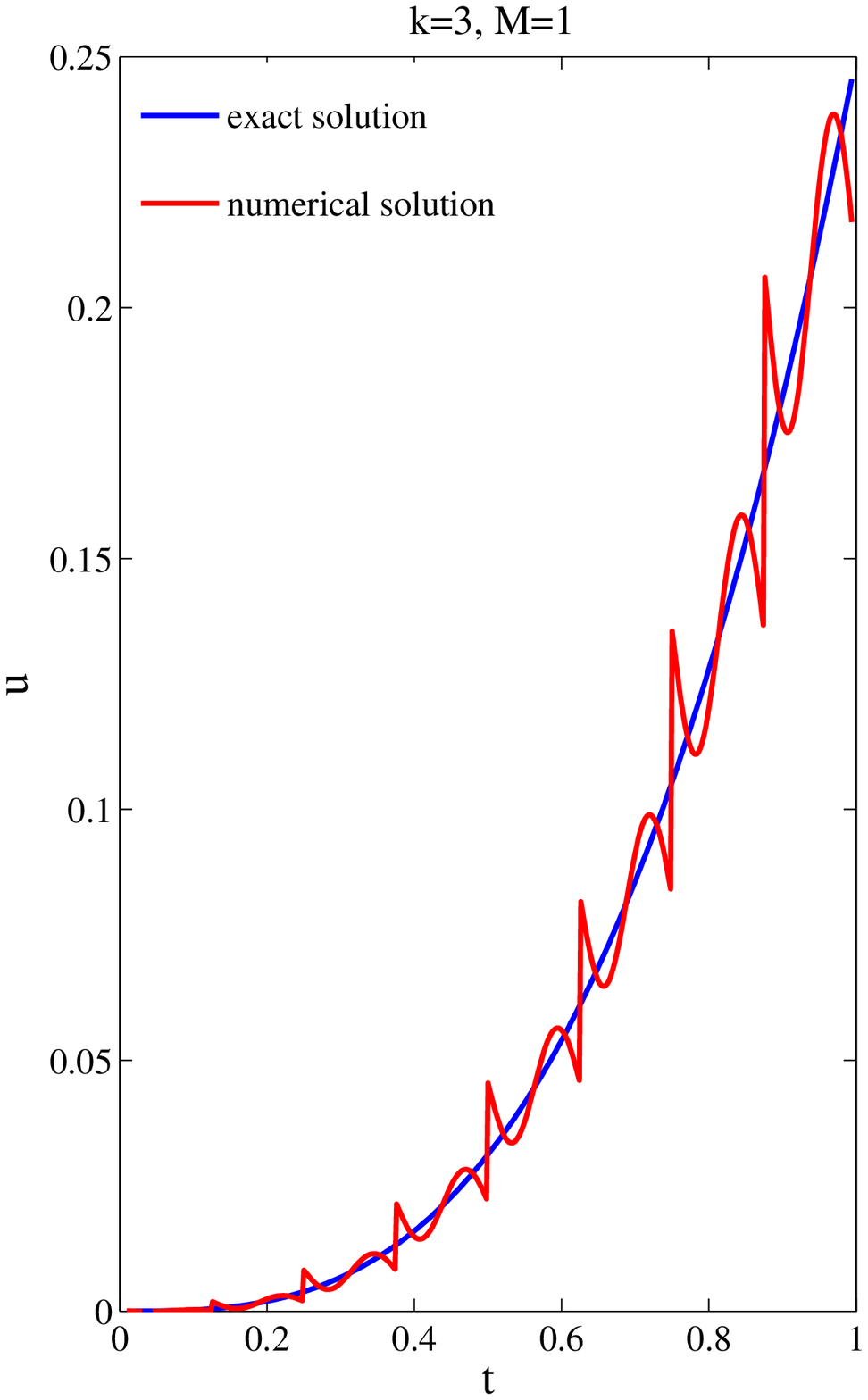}
\end{minipage}
\begin{minipage}[t]{0.32\linewidth}
\includegraphics[width=1.65in]{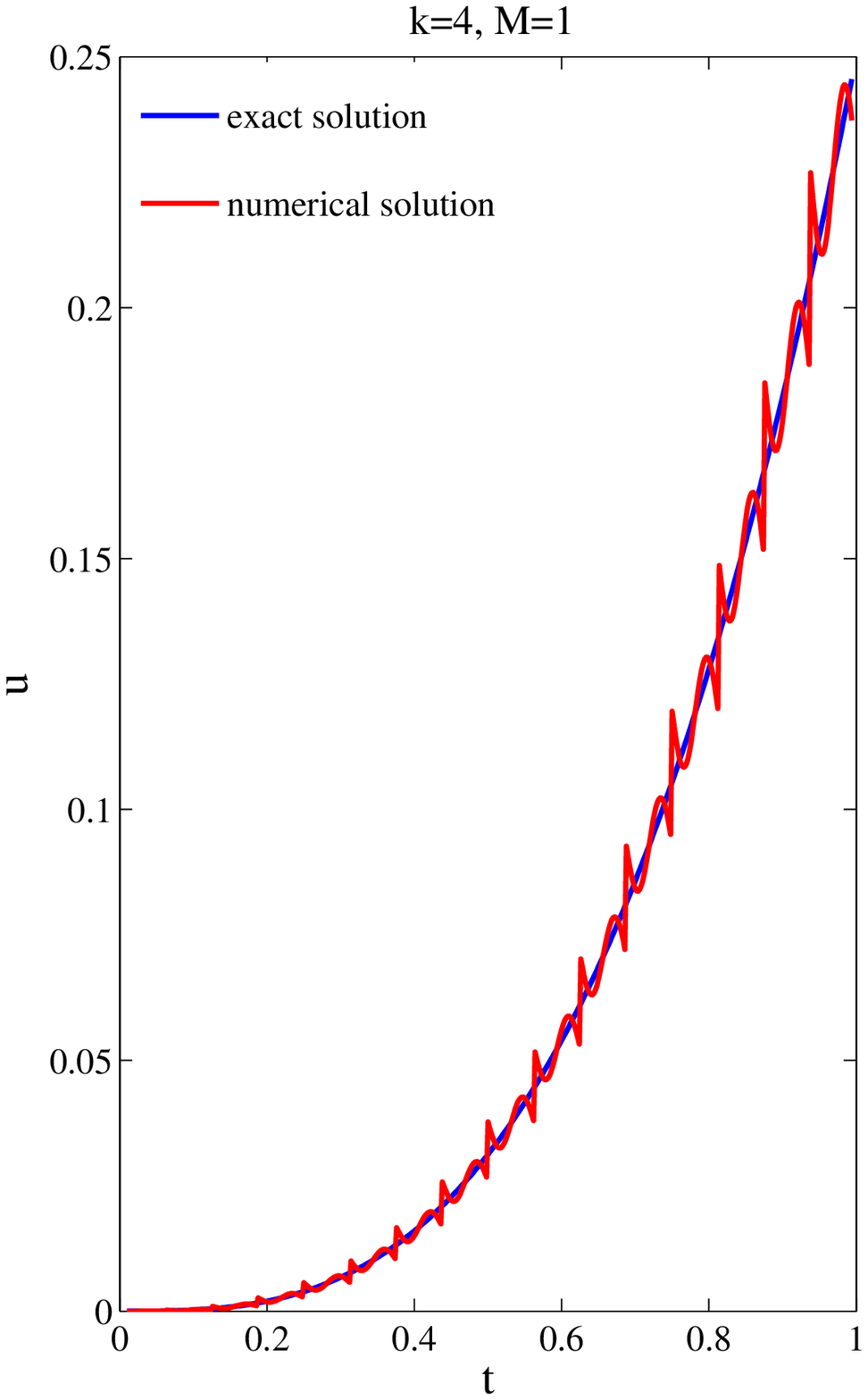}
\end{minipage}
\caption{The comparison between exact and approximate solutions when $\alpha=0.2$, $x=0.5$.}\label{fig1}
\end{figure}

\noindent
\textit{Example 6.1}
Let $a(x)=1$, $b(x)=-1$, $\ell=1$, and the initial-boundary values
$y(x,0)=0$, $y(0,t)=0$, $y(\ell,t)=t^3$. The right-hand function is selected as
\begin{equation*}
f(x,t)=\frac{6x^2t^{3-\alpha}}{\Gamma(4-\alpha)}-2t^3(1-x),
\end{equation*}
to give the exact solution $y(x,t)=t^3x^2$. Taking $\alpha=0.2$, $p=1$, and $N_h=20$,
Fig. \ref{fig1} and Table \ref{tab1} show the behavior of the approximate solutions versus the variation of $t$ at $x=0.5$
contrasted to the exact solutions and the absolute errors on several nodal points at $t=0.25$ for various $k$, $M$, respectively.\\

\begin{table*}%[!htb]
\centering
\caption{The absolute errors with $\alpha=0.2$ and various $k$, $M$ at $t=0.25$ for Example 6.1} \label{tab1}
\begin{tabular}{ccccc}
\toprule
x &$k=2$, $M=1$  &$k=3$, $M=1$ &$k=4$, $M=1$  &$k=5$, $M=1$  \\
\midrule0.1     &4.7632e-04  &1.7111e-04  &  7.1433e-05  &3.2587e-05    \\
        0.2     &1.9030e-03  &6.8340e-04  &  2.8502e-04  &1.2976e-04    \\
        0.3     &4.2801e-03  &1.5368e-03  &  6.4075e-04  &2.9153e-04    \\
        0.4     &7.6074e-03  &2.7315e-03  &  1.1386e-03  &5.1787e-04    \\
        0.5     &1.1885e-02  &4.2672e-03  &  1.7787e-03  &8.0880e-04    \\
        0.6     &1.7113e-02  &6.1441e-03  &  2.5608e-03  &1.1643e-03    \\
        0.7     &2.3291e-02  &8.3621e-03  &  3.4851e-03  &1.5844e-03    \\
        0.8     &3.0419e-02  &1.0921e-02  &  4.5515e-03  &2.0690e-03    \\
        0.9     &3.8498e-02  &1.3821e-02  &  5.7600e-03  &2.6182e-03    \\
\bottomrule
\end{tabular}
\end{table*}

\begin{figure}
\begin{minipage}[t]{0.49\linewidth}
\includegraphics[width=2.3in]{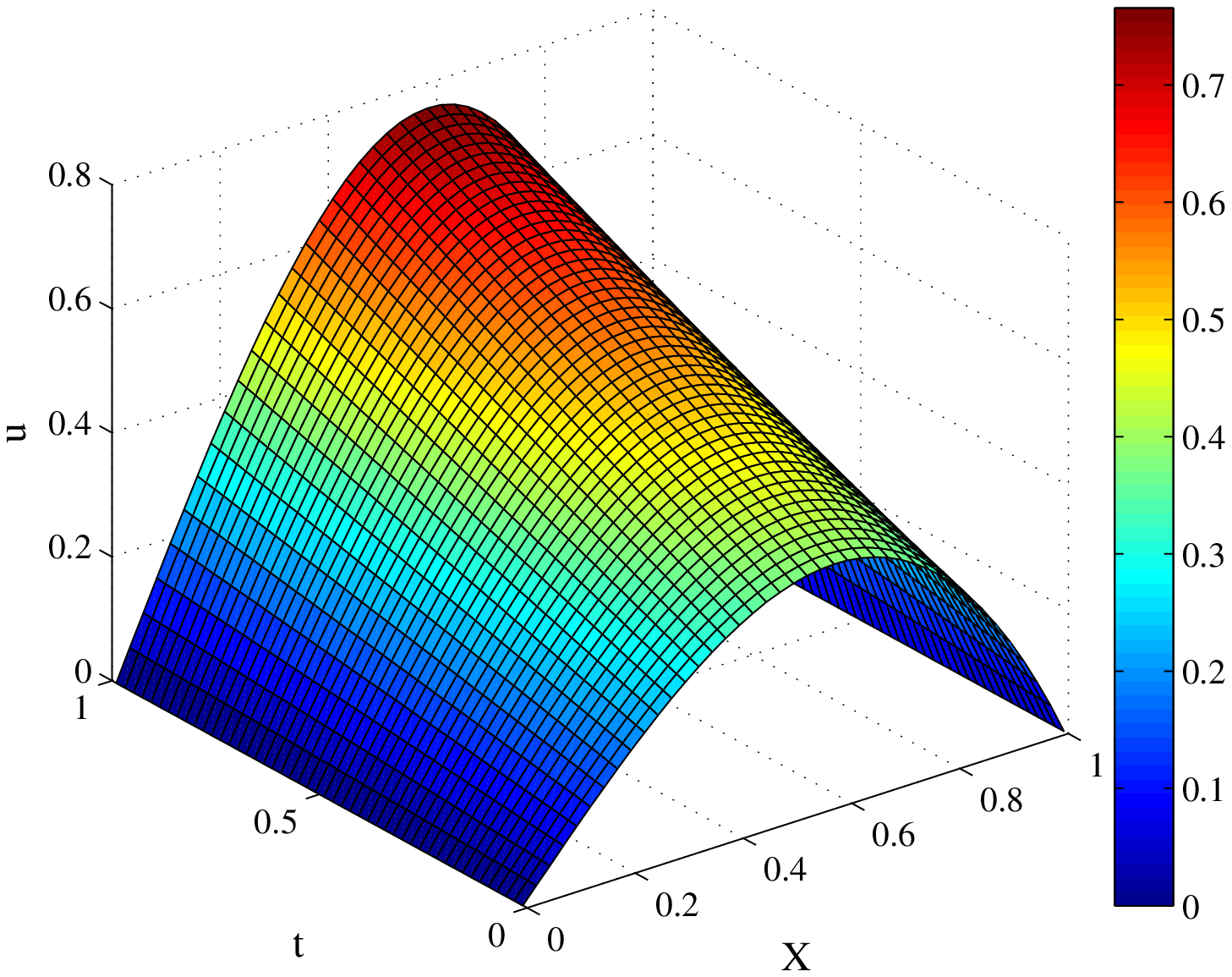}
\end{minipage}
\begin{minipage}[t]{0.5\linewidth}
\includegraphics[width=2.3in]{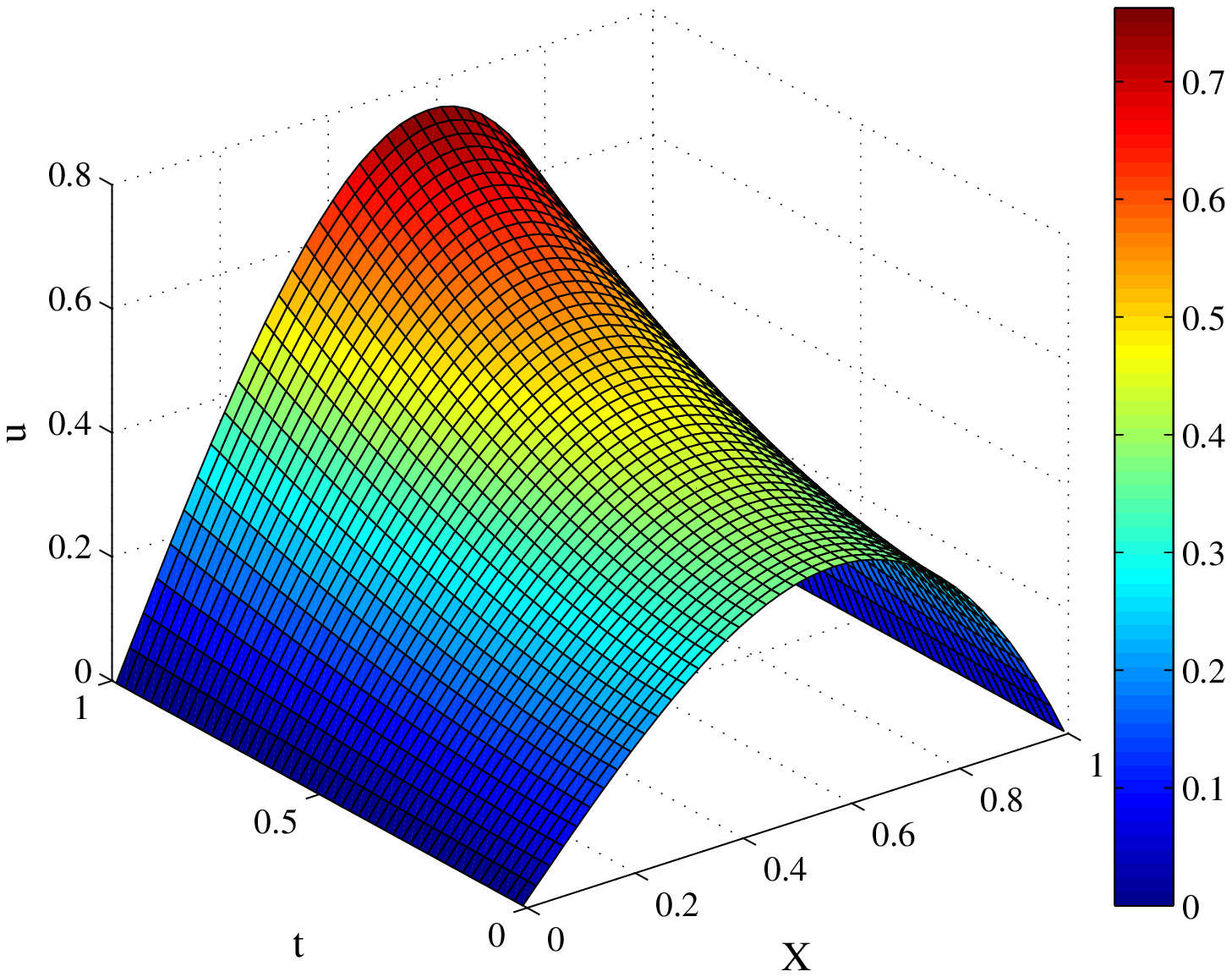}
\end{minipage}
\caption{The numerical solutions with $k=4$, $M=1$, $N_h=40$ for $\alpha=0.5$ and $\alpha=0.9$.}\label{fig2}
\end{figure}

\noindent
\textit{Example 6.2}
Let $a(x)=x$, $b(x)=-1$, $\ell=1$, and the initial-boundary values $y(x,0)=x-x^3$,
$y(0,t)=y(\ell,t)=0$. The forcing term and exact solution are %right side
\begin{equation*}
f(x,t)=\frac{\Gamma(1+2\alpha)}{\Gamma(1+\alpha)}t^\alpha(x-x^3)+(1+t^{2\alpha})(7x-3x^3),
\end{equation*}
and $y(x,t)=(1+t^{2\alpha})(x-x^3)$, respectively. The code is run with $p=1$. Fig. \ref{fig2} gives the numerical solutions
for $\alpha=0.5$, $\alpha=0.9$ by using $k=4$, $M=1$, and $N_h=40$.
Table \ref{tab2} reports the absolute errors at $t=0.5$ for various $k$, $M$ with $\alpha=0.7$ and $N_h=10$.
It is seen that our method is well convergent. \\

\begin{table*}%[!htb]
\centering
\caption{The absolute errors with $\alpha=0.7$ and various $k$, $M$ at $t=0.5$ for Example 6.2} \label{tab2}
\begin{tabular}{ccccc}
\toprule
x &$k=4$, $M=1$  &$k=5$, $M=2$ &$k=6$, $M=1$  &$k=7$, $M=2$  \\
\midrule0.1     &3.3203e-03  &1.6340e-03  &  8.0029e-04  &3.8732e-04    \\
        0.2     &6.4390e-03  &3.1686e-03  &  1.5517e-03  &7.5082e-04    \\
        0.3     &9.1546e-03  &4.5045e-03  &  2.2055e-03  &1.0667e-03    \\
        0.4     &1.1266e-02  &5.5425e-03  &  2.7130e-03  &1.3114e-03    \\
        0.5     &1.2571e-02  &6.1836e-03  &  3.0257e-03  &1.4615e-03    \\
        0.6     &1.2870e-02  &6.3291e-03  &  3.0955e-03  &1.4938e-03    \\
        0.7     &1.1961e-02  &5.8808e-03  &  2.8746e-03  &1.3855e-03    \\
        0.8     &9.6461e-03  &4.7409e-03  &  2.3158e-03  &1.1146e-03    \\
        0.9     &5.7250e-03  &2.8126e-03  &  1.3728e-03  &6.5958e-04    \\
\bottomrule
\end{tabular}
\end{table*}

\begin{figure}
\begin{minipage}[t]{0.49\linewidth}
\includegraphics[width=2.3in]{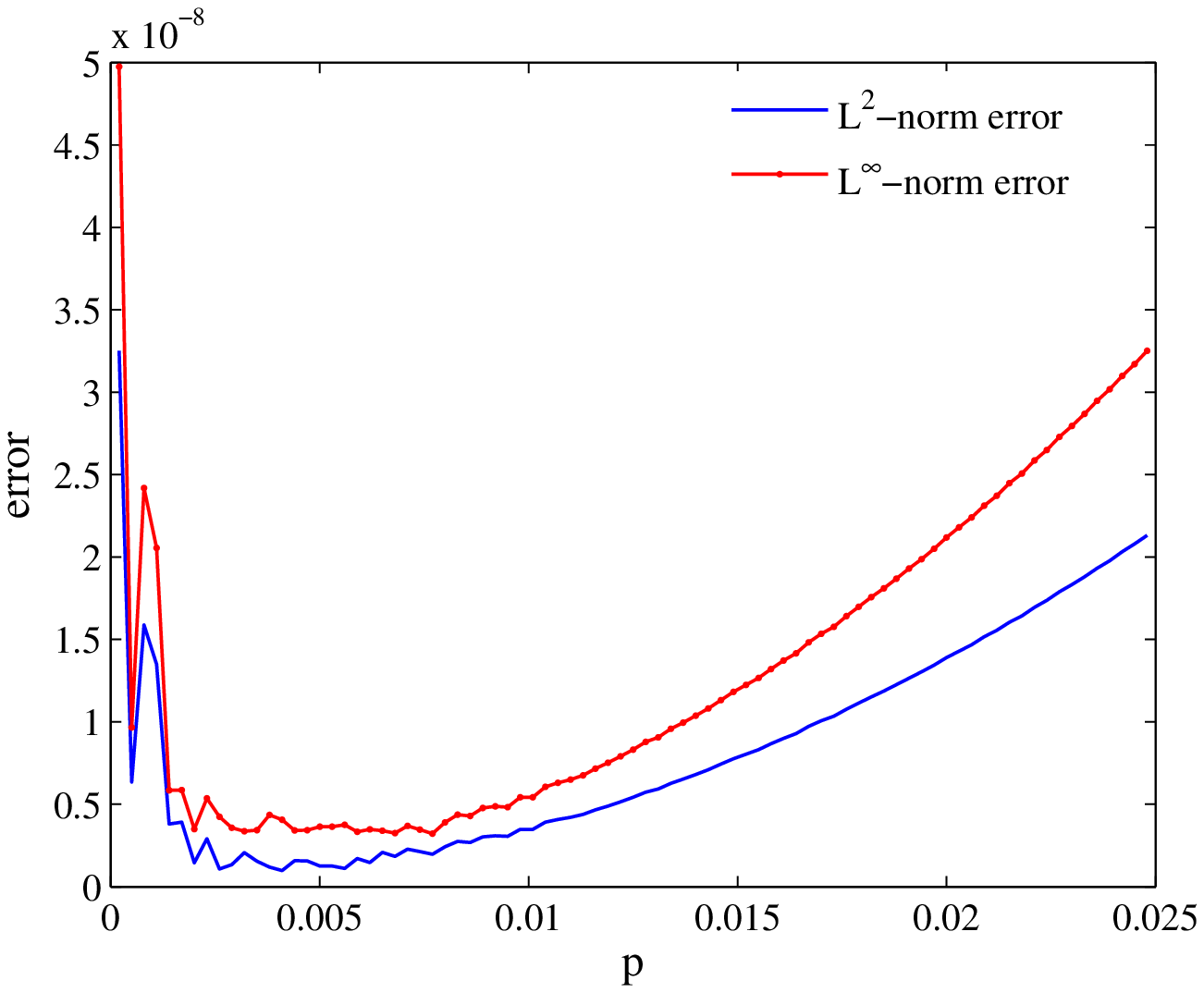}
\end{minipage}
\begin{minipage}[t]{0.5\linewidth}
\includegraphics[width=2.3in]{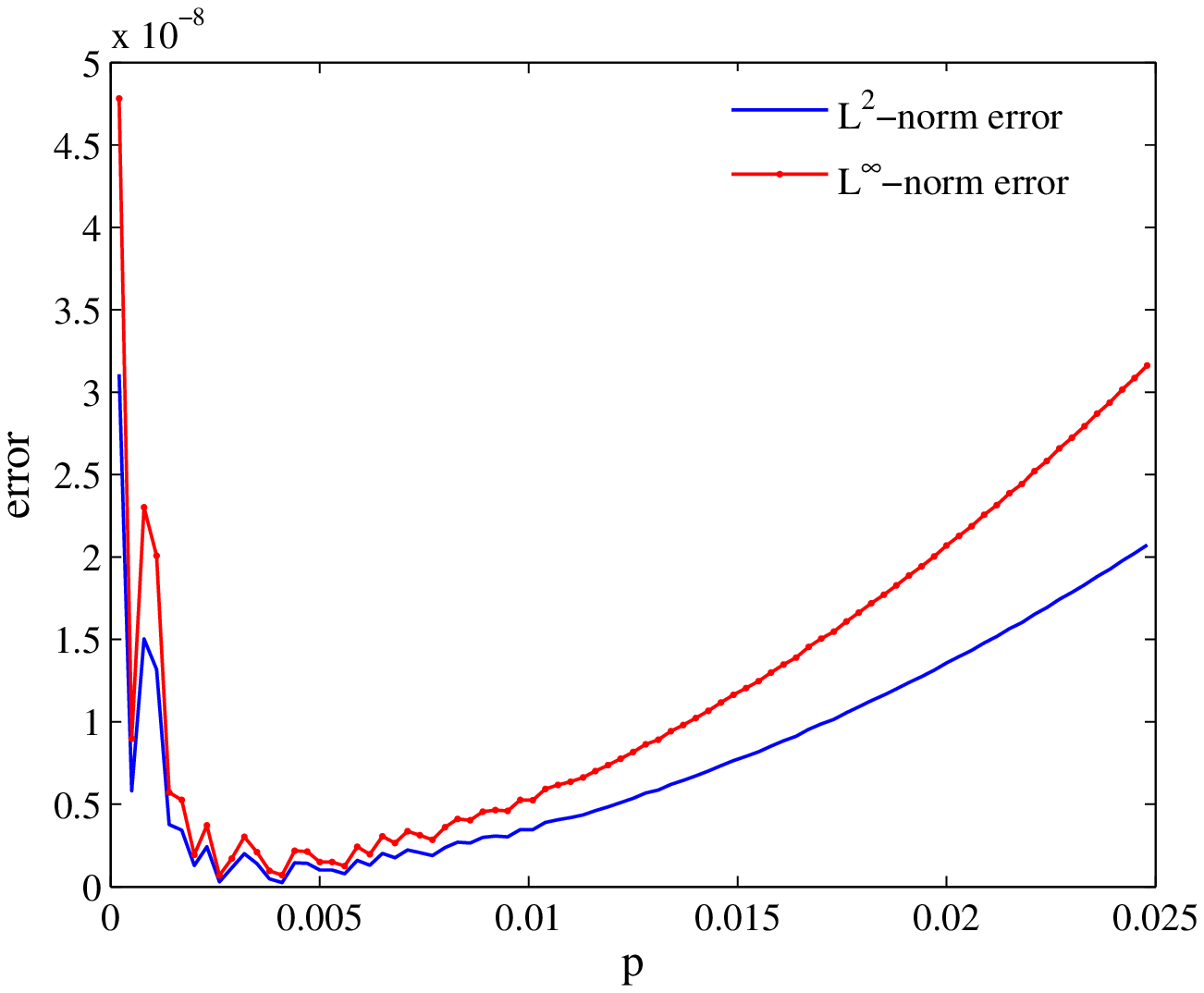}
\end{minipage}
\caption{The global errors versus the variation of $p$ with $k=6$, $M=1$, $N_h=20$ for $\alpha=0.3$ and $\alpha=0.6$.}\label{fig3}
\end{figure}

\begin{figure}
\begin{minipage}[t]{0.49\linewidth}
\includegraphics[width=2.3in]{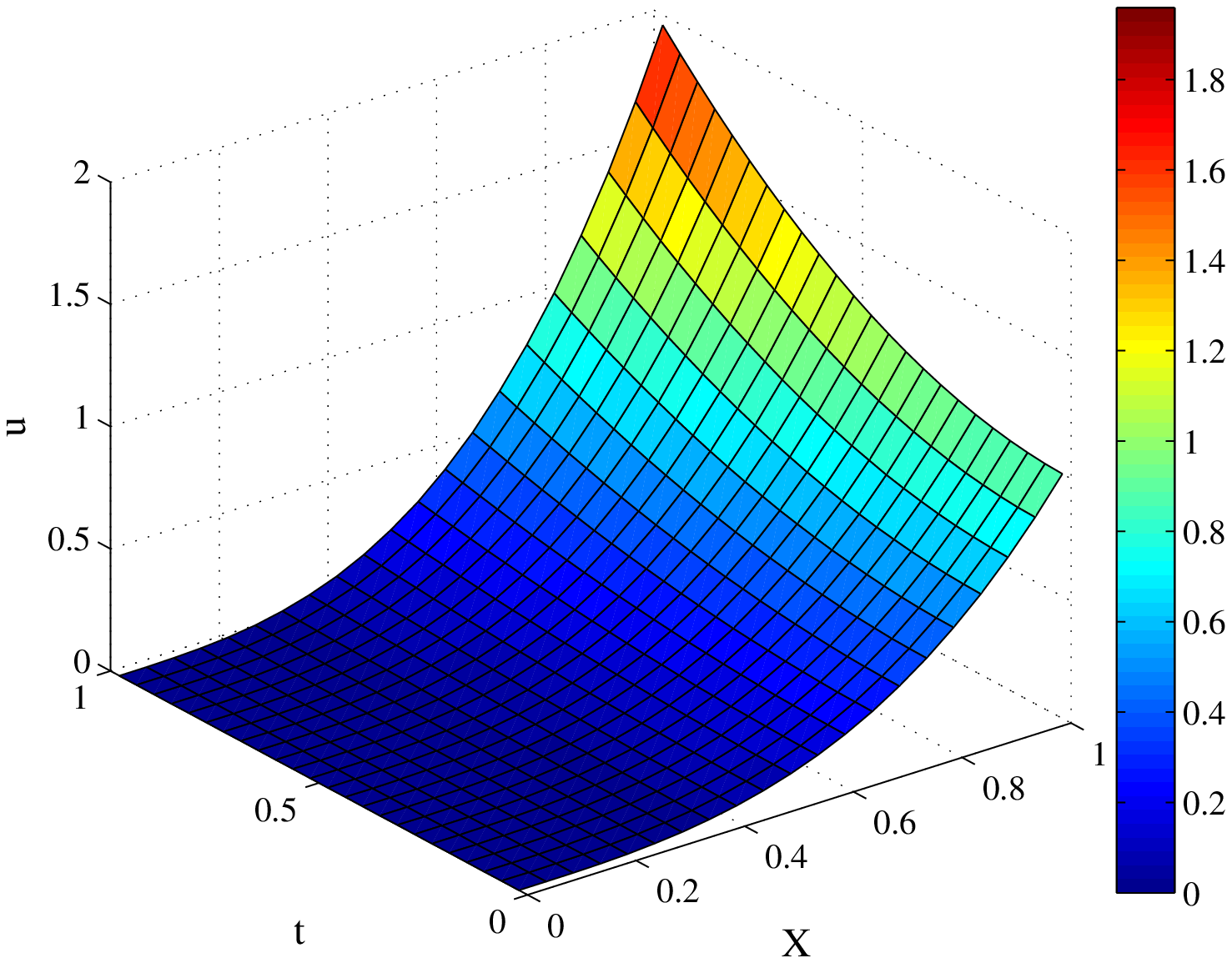}
\end{minipage}
\begin{minipage}[t]{0.5\linewidth}
\includegraphics[width=2.3in]{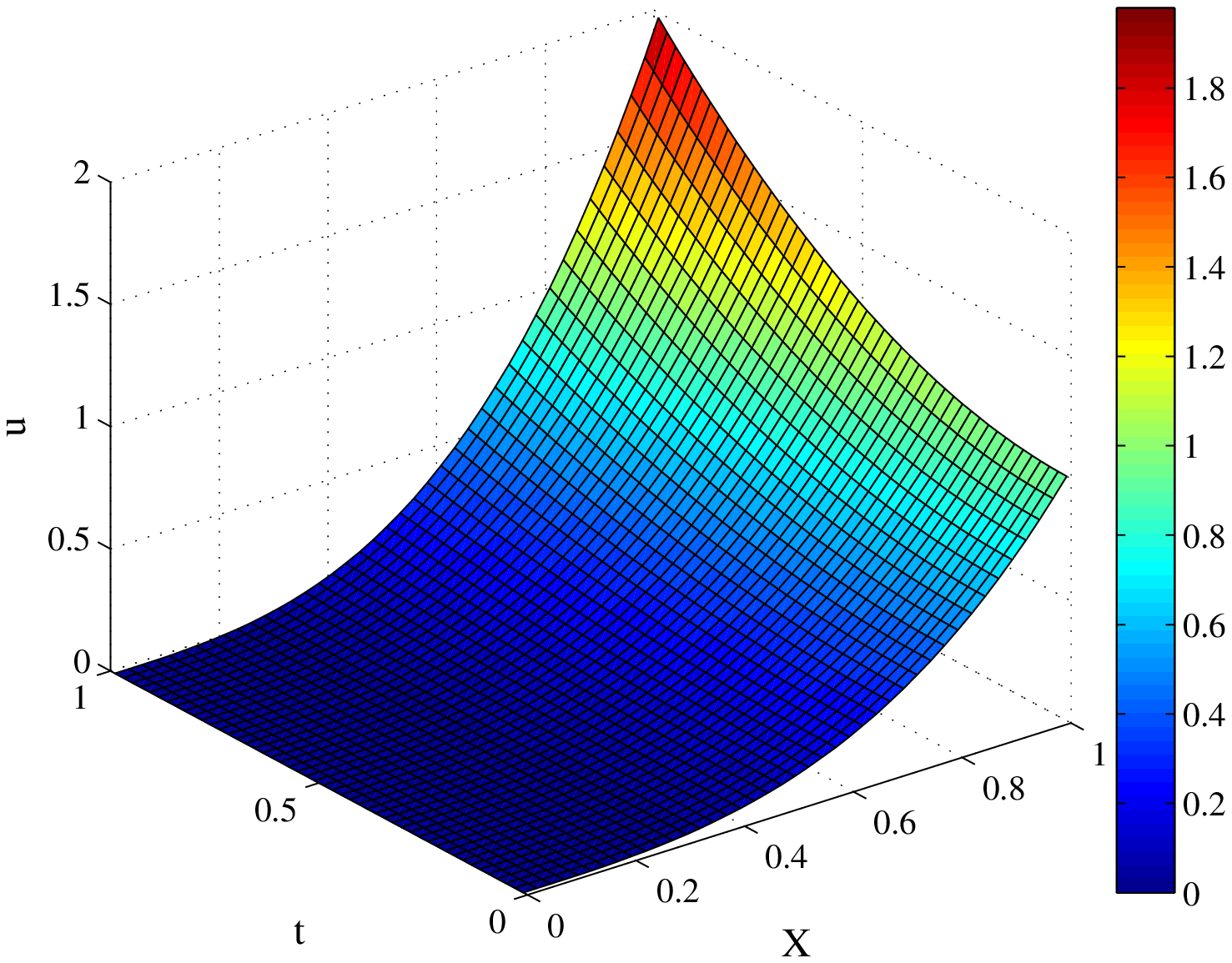}
\end{minipage}
\caption{The numerical solutions with $k=3$, $M=1$, $N_h=20$ and $k=4$, $M=1$, $N_h=40$ for $\alpha=0.1$.}\label{fig4}
\end{figure}

\noindent
\textit{Example 6.3}
In this test, we consider Eqs. (\ref{eq01})-(\ref{eq03}) with
$a(x)=1$, $b(x)=-x$, $\ell=1$, and the initial-boundary values $y(x,0)=x^3$,
$y(0,t)=0$, $y(\ell,t)=1+t^2$, and the forcing function
\begin{equation*}
   f(x,t)=\frac{2t^{2-\alpha}x^3}{\Gamma(3-\alpha)}-3(1+t^2)x^2.
\end{equation*}
The exact solution is given by $y(x,t)=(1+t^2)x^3$. In Fig. \ref{fig3}, we depict the global errors versus the
variation of parameter $p$ for $\alpha=0.3$, $\alpha=0.6$ with $k=6$, $M=1$, and $N_h=20$.
As the graph show, the optimal $p$ for this problem is roughly located on $[0.002,0.005]$.
On taking $p=0.025$, the numerical solutions for $\alpha=0.1$ computed by using $k=3$, $M=1$, $N_h=20$
and $k=4$, $M=1$, $N_h=40$ are displayed in Fig. \ref{fig4}, and some numerical results at $t=t_{N_t}$ with $k=6$, $M=2$
are tabulated in Table \ref{tab3}, where the second-order convergent rates and good accuracy are observed and confirmed. \\

\begin{table*}[!htb]
\centering
\caption{The numerical results with $k=6$, $M=2$ and various $\alpha$ at $t=t_{N_t}$ for Example 6.3.} \label{tab3}
\begin{tabular}{cccccc}%{llllll}
\toprule
          $\alpha$  &  $N_h$   &$e_2(t_{N_t},N_h)$ &Cov. rate &$e_\infty(t_{N_t},N_h)$ &Cov. rate  \\
\midrule $0.3$      & 5        &3.4402e-07  &-       &    5.0737e-07  &-     \\
                    & 10       &8.6712e-08  &1.9882  &    1.3205e-07  &1.9419   \\
                    & 20       &2.1705e-08  &1.9982  &    3.3023e-08  &1.9996   \\
                    & 40       &5.3797e-09  &2.0125  &    8.2022e-09  &2.0094   \\
       $0.6$        & 5        &3.3183e-07  &-       &    4.8881e-07  &-     \\
                    & 10       &8.3584e-08  &1.9891  &    1.2771e-07  &1.9364   \\
                    & 20       &2.0960e-08  &1.9956  &    3.2003e-08  &1.9966   \\
                    & 40       &5.2376e-09  &2.0006  &    8.0132e-09  &1.9978   \\
       $0.9 $       & 5        &3.1876e-07  &-       &    4.6891e-07  &-     \\
                    & 10       &8.0489e-08  &1.9856  &    1.2333e-07  &1.9268   \\
                    & 20       &2.0495e-08  &1.9735  &    3.1339e-08  &1.9765   \\
                    & 40       &5.4419e-09  &1.9131  &    8.3056e-09  &1.9158  \\
\bottomrule
\end{tabular}
\end{table*}

\section{Conclusion}
In this article, we have presented a numerical algorithm for the variable coefficient
time-fractional convection-diffusion equation by approximating its solution as a truncated sine-cosine wavelets series
in time and exponential spline interpolation in space. %together with powerful  utilized treat
The method turns a differential model into a solvable algebraic system via
the derived wavelet operational matrix of fractional integration and thereby is of significance to
make the execution of program easier and more economic since the historical correlation of these models.
The tested results of its codes on some illustrative examples have manifested that it can handle %solves
the given equations very effectively. It is certain that our method
%can serve as an alternative choice to simulate the other fractional problems.
provides an alternative to simulate the other fractional problems.

\begin{acknowledgements}
This research was supported by National Natural Science Foundations of China (No.11471262 and 11501450).
\end{acknowledgements}
%\small
%\textbf{Conflict of interest} The authors declare that this work is our original research and  there is no conflict of interests.

\bibliographystyle{spmpsci}      % mathematics and physical sciences
\bibliography{reference}   % name your BibTeX data base

\end{document}